\input amstex
\documentstyle{amsppt}
\NoRunningHeads 
\magnification = 1000   
\NoBlackBoxes      
\topmatter
\title  Indecomposable Higher Chow Cycles 
on low dimensional Jacobians
           \endtitle

\author {Alberto  Collino}   \endauthor 

\affil { \it Dipartimento di Matematica,
Universit\`a di Torino, \\ Via Carlo Alberto 10,
10123 Torino, Italy \\
e-mail: collino\@dm.unito.it  }\endaffil
\thanks Research partially supported  European Research
Project AGE-Algebraic Geometry in Europe  and by   Progetto di ricerca di
interesse nazionale MURST-1997 su Geometria algebrica, Algebra commutativa
e Aspetti computazionali  \endthanks 
\abstract\nofrills\ 
There is a basic indecomposable higher cycle $K \in CH^{g}(J(C),1)$ 
on the Jacobian $J(C)$ of a general
hyperelliptic curve $C$ of genus $g$, see \cite{3}. 
Consider $K_t$  the translation of $K$ associated with a point $t \in C$,
we prove that in general $K - K_t$ is indecomposable if $g \geq 3$. 
Our tool is Lewis' condition   
for indecomposability \cite{8}.  We also show that on the jacobian $J(C)$ of 
any curve $C$ 
of genus $3$ there is a  geometrically natural family of higher cycles,
when $C$ becomes hyperelliptic the family in the  limit contains  a component of 
indecomposable cycles of type $K - K_t$.
\endabstract
\endtopmatter

\document 
\head Introduction
\endhead

	The first higher Chow group $ CH^{p}(X,1) \simeq H^{p-1}(X,\Cal K_{p})$ 
of a non singular variety $X$ is generated by higher cycles of the form
$ Z = \sum_i Z_i \otimes f_i$, where the  $Z_i$ are irreducible
subvarieties of codimension $(p-1)$ and the rational functions
$f_i \in \Bbb C (Z_i)^\times$ obey the rule $\sum_i div(f_i) =0$ as a cycle on  $X$.  
		Consider the subgroup of {{\it} decomposable\/} cycles 
$$
CH_{\text{dec}}^k(X,1) := Im \big\{ CH^1(X,1) {\otimes}  CH^{k-1}(X)
\longrightarrow CH^k(X,1) \big\} \quad,
$$
and the related quotient of {{\it} indecomposable\/} cycles 
$$
CH_{\text{ind}}^k(X,1;{\Bbb Q}) := 
\big( CH^k(X,1)\big/ CH_{\text{dec}}^k(X,1)\big) {\otimes}
{\Bbb Q} \quad.
$$
	A geometrically natural higher cycle $K$ for $CH^g(J(C),1)$ 
was given in \cite{3}, here  $J(C)$ is a hyperelliptic Jacobian, and it was proved
that $K$ is indecomposable  when $C$ is general. This result was reached 
by showing that 
the interesting part of the regulator image of $K$ must be non trivial. 
The regulator or cycle class is a map 
$cl_{k,1}$  from
the higher Chow groups into Deligne cohomology \cite{1}, 
$$
 cl_{k,1} :CH^k(X,1;\Bbb Q)  \longrightarrow H_{\Cal D}^{2k-1}(X,\Bbb Q (k)) \, , 
$$
\noindent which satisfies a  rigidity property of Beilinson type \cite{10}, 
namely $ CH_{\text {ind}}^k(X,1;\Bbb Q)$ has countable image
modulo `Hodge classes'. 
	Our aim here is to show that  
translation on $J$ acts non trivially on $K$ when the genus of $C$ is (at least)
$3$, to the effect that the orbit of $K$ in  $CH_{\text{ind}}^g(J(C),1;{\Bbb Q})$ 
is uncountable. We use in an essential way Lewis' condition 
for indecomposability \cite{8}.
	Voisin has conjectured some time ago that for a smooth projective 
variety $X$ the group $CH_{\text{ind}}^2(X,1;\Bbb Q)$ is countable \cite{11}, 
we verify below  that indeed in genus $2$ translation does not change the 
indecomposable class of $K$.

	Lewis' conditions are some Hodge-theoretic properties of the regulator image
of a complete family  of higher cycles on $X$, when they hold
the general member in the family  is indecomposable, even if its regulator class
vanishes. The main result of \cite{6} shows that indeed there is a complete family
satisfying the conditions, in their example $X$ is a sufficiently 
general product of three elliptic curves.	The cycles involved may appear to be ad-hoc, 
this has moved me to consider certain
geometrically natural cycles, the '4-configurations', 
which live on the general jacobian of genus $3$.
By specializing to the hyperelliptic locus a 4-configuration becomes 
the hyperelliptic configuration $K - K_t$, that is the difference of the basic cycle 
by its translation associated  with a point $t \in C$. 
This leads us to  conjecture that in general the 4-configuration is indecomposable. 

	In the  last part we show that translation on a bielliptic jacobian $J(G)$
acts non trivially on the indecomposable type of some cycles in  $CH^{3}(J(G),2)$.
This is in contrast with the known fact that on an elliptic curve $E$ translation acts
trivially on  $ CH^{2}(E,2) \simeq H^{0}(E,\Cal K_{2})$  modulo $K_{2}(\Bbb C)$.

\heading 1. Lewis' condition. \endheading

	In this section we recall some notation and definitions and then
we state one of  Lewis'  theorems. Our aim
is to have a concrete reference at hand, 
for more details the  reader should consult either the original paper \cite{8} 
or the excellent survey \cite{7}.

\subhead 1.1 The real regulator and some definitions
\endsubhead

	Let   $X$ be projective and nonsingular, it is 
$$
\frac{H^{i-1}(X,\Bbb C)}{F^jH^{i-1}(X,\Bbb C) + H^{i-1}(X,\Bbb R(j))} \simeq
\frac{H^{i-1}(X,\Bbb R(j-1))}{\pi_{j-1}(F^jH^{i-1}(X,\Bbb C))},
$$
and therefore  one has
$$\align
 H_{\Cal D}^{2k-1}(X,\Bbb R(k)) &\simeq H^{2k-2}(X,\Bbb R(k-1)) \cap
F^{k-1}H^{2k-2}(X,\Bbb C) \\ \vspace{1\jot}
  &=: H^{k-1,k-1}(X,\Bbb R(k-1)) .
\endalign
$$

According to  Beilinson \cite{1}  the  real  regulator image  of a cycle
$ Z \in CH^k(X,1;\Bbb Q)$  is the element
$$
 R_{k,1}( Z) \in  H_{\Cal D}^{2k-1}(X,\Bbb R(k)) \simeq H^{k-1,k-1}(X,\Bbb R(k-1)),
$$
determined by the class of the current
$$
 R_{k,1}( Z) : \omega \longmapsto (2\pi \sqrt{-1})^{k-1-d} \sum_i \int_{Z_i -
Z_i^{\text{sing}}} \omega \,\log|f_i| \ .
$$

\definition {  Definition:   level of a quotient $ CH^k(X,m;\Bbb Q)\big/ G$}

Consider  a subgroup  $G$  of $ CH^k(X,m;\Bbb Q)$,
the { \it level } of $ CH^k(X,m;\Bbb Q)\big/ G$ is the
least nonnegative integer ~$r$ such that
$$
  CH^k(X,m;\Bbb Q) = G + i_* CH^{r+m}(Y,m;\Bbb Q)
$$
 for some $i:Y \hookrightarrow X$ closed and of pure codimension ~$(k-r-m)$.
\enddefinition

\definition{  Definition:  $ H_N^{k-l,k-m}(X)$}

 	The {\it  coniveau filtration } is given by
$$
 N^jH^i(X,\Bbb Q) := \ker \big(\,H^i(X,\Bbb Q) \longrightarrow lim_{codim{_X} Y\ge j}
H^i(X-Y,\Bbb Q) \,\big),
$$
 where the direct limit is over closed subvarieties $Y \subset X$.  We
let
$$H_N^{k-l,k-m}(X) :=Im (\, N^{k-l}H^{2k-l-m}
(X,\Bbb Q)\otimes \Bbb C \longrightarrow
H^{k-l,k-m}(X) \,\big) \quad,$$

and we note that it is the complex subspace generated by 
the Hodge projected image of the
coniveau filtration.
\enddefinition

\remark{Fact}  
	Lewis constructs certain  complex subspaces 
$$H^{\{k,l,m\}}(X) \subseteq  H^{k-l,k-m}(X) \quad ,$$
such that for $m=0$  one has 
$$H^{\{k,l,0\}}(X) \subseteq  H_N^{k-l,k}(X)$$ 
\endremark
	The spaces $H^{\{k,l,m\}}(X)$ are obtained by a process of
K\"unneth projection of the Hodge components 
of the real regulator classes of the elements in $CH^k(X \times S;m)$,
with varying $S$, $S$ projective and nonsingular.

 \subhead 1.2 Lewis's theorem
\endsubhead
We recall  Lewis's theorem in an abridged version,
because we need it only for  $CH_{\text{ind}}^k(X,1;{\Bbb Q})$.
\proclaim{ Theorem {\rm (\cite{8})}}
 Let $X$ be a projective algebraic manifold. Then
\roster
\item $H^{\{k-1,l-1,0\}}(X) \subset H^{\{k,l,1\}}(X)$.
\smallskip
\item If $H^{\{k,l,1\}}(X) \big/ H^{\{k-1,l-1,0\}} (X) \ne 0$ then
${  level}( CH_{\text{ind}}^k(X,1;\Bbb Q)) \ge l-1$.
\endroster 
\endproclaim

\proclaim{ Corollary}
 If $H^{\{k,l,1\}}(X) / H^{\{k-1,l-1,0\}} (X) \ne 0$ for some  $l$
with $2 \le l \le k$, then $CH_{\text{ind}}^k(X,1;\Bbb Q)$ is uncountable.
\endproclaim

We shall prove that the hypotheses of the corollary
hold when $X$ is the general hyperelliptic jacobian of genus $3$ for $k=3$, $l=2$. 
 
\heading 2. The basic hyperelliptic construction. \endheading
\subhead 2.1. The hyperelliptic configuration
\endsubhead
Let  $f : C  \to \Bbb P^{1}$   
be the double cover associated with a 
hyperelliptic curve, we fix 
two ramification points $w_1$ and $w_2$  
on $C$ and choose a standard parameterization   
on $\Bbb P^{1}$  so that $ f(w_1)=0$
and $ f(w_2) = \infty $. 
The points $w_1$ and $w_2$ are for us the distinguished Weierstrass points,
and $ \epsilon  := class({w_1}-{w_2})$
is the associated element of order two in $Pic^0(C) $.
 
	It is  convenient  for us to identify  $J(C)=Pic^0(C) $ with $Pic^1(C)$ by adding 
$w_1$. We embed $C$  in the natural way  in  $Pic^1(C)$, for $t \in Pic^1(C)$  
we let $C_t$ be the translate of $ C $ by class$(t-w_1)$. 
Consider  now $ W_1 := C = C_{w_1}$ and $ W_2 :=C_{w_2}$,
the $ \epsilon $ translate of $C$, and fix a point $t \in C$. 
It is useful to observe that
the intersection  $C_t \cap W_1 \cap W_2$
is the point  $w_1$. 
We shall follow the convention to indicate a rational function on $C_t$
by using the same name given to the corresponding function on $C$.
  
	Consider $K :=   W_1 \otimes f  +  W_2 \otimes f $ and  its 
t-translation $K_t :=   C_t \otimes f  +  C_{t+ \epsilon}\otimes f $. Now
$K$ is the basic cycle of \cite {3}, hence it is a non trivial indecomposable 
element of $CH^{g}(J(C),1;\Bbb Q)$ for $C$ general enough.
Here we are interested in  the cycle $Z_t := K - K_t$, which we call 
the hyperelliptic configuration. 
It is clear that $Z_t$ belongs to the kernel of the regulator map,
because the regulator at hand is computed as a current acting on 2-forms 
on $J$, hence it is indifferent to translation. 
	Our aim is to prove that for general $t$ the hyperelliptic configuration
$Z_t $ is  indecomposable if the genus is at least $3$, but first we show
that on the contrary 
\proclaim{Proposition} In genus $2$ the hyperelliptic configuration is decomposable.
\endproclaim 
 
The divisor class of $W_1-W_2$ is of order $2$ in $Pic^0(J(C))$ hence there is
a rational function $\theta$, say, with $div(\theta)$ $= $  $2W_1-2 W_2$, and
we may assume without restriction that  $\theta$ is symmetric with respect to the $\epsilon$
involution on $J(C)$, namely that  $\epsilon^* \theta = \theta^{-1}$.
The other two curves in the support of the configuration $C_t$ and $C_{t+\epsilon}$
are  the translation respectively of $W_1$ and $W_2$  
and thus $2 (C_t- C_{t+\epsilon})$ $= $  $div(\theta_t)$, where $\theta_t$ is also
symmetric. The divisor of $(\theta_t)_{|W_ 1}$ 
is supported on $W_1$ at the distinguished  Weierstrass points because at the other
point $t$ of intersection of $W_1$ with $C_t$ the contributions cancel out,
and in fact $(\theta_t)_{|W_ 1}$ $=$ $c_tf$ 
for some nonzero constant $c_t$, depending on $t$,
because we see that $(\theta_t)_{|W_ 1}$ and $f$ have the same divisor.
 
By symmetry  $(\theta_t)_{|W_2} =  (c{_t} f)^{-1}$,
and also
$(\theta)_{|C_t}$$=$ $k_tf$
and $(\theta)_{|C_{t+\epsilon}}$$=$ $(k{_t} f)^{-1}$.
	By definition the tame symbol
of  two rational functions $a$ and $b$   yields  at a divisor $D$  the rational function
$ (T \{a,b \})_D = (-1)^{{\nu_{D}(a)} {\nu_{D}(b)}}
(a^{ \nu_{D}(b)} /b^ {\nu_{D}(a)} )_{|D} $
where $\nu_{D}$ is the order at $D$.
	In this way it is
$$ T \{\theta_t,\theta \}  
 =   2 (W_1 \otimes f   +  W_2 \otimes f
+  C_t\otimes f ^{-1}  +  C_{t+{\epsilon}} \otimes f ^{-1}  + 
(W_1+  W_2) \otimes c_t - (C_t  + C_{t+{\epsilon}}) \otimes k_t)  ,$$
and therefore $Z_t$  vanishes in 
$CH_{\text{ind}}^2(J(C),1;{\Bbb Q})$.

\subhead 2.2. A higher Chow cycle on $J\times C$
\endsubhead
We keep the identification $J(C) \simeq Pic^1(C)$.
We  embed  $C\times C$ in $J\times C$ in four different ways,
two of them twisted  and two of them straight. 
The straight surfaces are   $S_1 :=  W_1\times C$ and $S_2 :=  W_2\times C$.
The twisted ones are $T_1$ and $T_2$, they
are the images in $J\times C$
of the maps which send a point $ (x,t)$ in $C\times C$ to
$(class(x  +t-w_1),t)$ and   $(class(x+t-w_2) ,t)$ respectively.
Consider the rational function
$F$ on $C\times C$  which is the pull-back under the first projection of $f$ on $C$.
We have on $S_i$ and $T_i$  rational functions $F_i$ and $G_i$ respectively,
they  correspond to $F$ under the said isomorphisms.

 In this manner ${\frak S } :=   S_1{\otimes} F_1 + S_2{\otimes} F_2$ and  
${\frak T} :=   T_1{\otimes} G_1 + T_2{\otimes} G_2$
are higher Chow cycles for $CH^{g}(J\times C,1;\Bbb Q)$.
Our final cycle is $${\frak Z } := {\frak S} -  {\frak T} \qquad  .$$ 
Taking sections over  $t \in C$ yields  
 ${\frak S }_{\ast}(t)  =  K$  and  ${\frak T }_{\ast}(t)  =  K_t$,
and thus we have ${\frak Z }_{\ast}(t)   = Z_t$,
the hyperelliptic configuration.

\subhead 2.3. Lewis' method and a theorem of Lieberman
\endsubhead

Lewis' conditions  deal with  the Hodge-K\"unneth decomposition
of the real regulator image of a higher cycle like ${\frak Z }$,
if they are satisfied then the general section  ${\frak Z }_{\ast}(t)$ is
indecomposable. We  prove below that this is the case
in our situation and  thus we have that but possibly for a countable subset of points in $C$: 
\proclaim{2.3.1. Theorem}
The hyperelliptic configuration  $ Z_t$ is indecomposable
\endproclaim

In order to proceed we need to recall that the General Hodge Conjecture asserts that
the filtration on cohomology with rational coefficients induced by the
Hodge filtration should coincide with the coniveau filtration.
Lieberman \cite{9} proved that the General Hodge Conjecture
holds for $H^3(J)$, when $J$ is a general hyperelliptic jacobian of genus $3$. 
More precisely the only
subhodge  structures in $H^3(J)$ are the primitive part  $P^3$
and the product  $\Theta H^1(J)$,  $\Theta$ being the class of the theta divisor.
In this way  the coniveau $1$ filtration $N^{1}H^3(J)$ is exactly $\Theta H^1(J)$ 
and therefore it is  orthogonal  to $P^3$.
	This is the point where we need to use the hypothesis $g=3$. 
Pirola (oral communication) claims that the analogue
of Lieberman's result holds in fact for any general hyperelliptic jacobian of genus at least 
$3$. If this is the case then our proof would yield 
indecomposability of $Z_t$ also for $g \geq 4$.

	In our situation   Lewis' condition  is fulfilled 
if we  prove that there exists a certain class $ \eta \in H^{2,1}(J) $ 
and a form $\nu \in H^{0,1}(C)$ such that  
the pairing   $ <R({\frak Z }), \eta \wedge \nu> \ne
0 $, where $R({\frak Z })$ is the real regulator image  of ~$ {\frak Z }$ 
defined as a current. In order to satisfy Lewis' requirement  $ \eta $ must be orthogonal
to $H^{ \{2,1,0\} }(J)$, and this is the case when $ \eta $ is primitive,
because of Lieberman's result.

\subhead 2.4. On the   real regulator image of $K$
\endsubhead

The standard inner product on the space $ H^0(C,K_C)$ of global holomorphic forms 
on a curve $C$ is $<\alpha,\beta> :=  (i/2) \int _ C \alpha \wedge \bar{\beta} $.
We fix a basis  $ \omega_i^{J} , \dots, \omega_g^{J}$ for the $1$-forms on $J$
with the property that   
$\zeta_i :=a^{*}\omega_i^{J}$ is an orthonormal basis for $ H^0(C,K_C)$,
here $a : C \to J$   is the standard map. 
The class of the theta divisor  on $J$ is then given by the form
$ \theta_J = (i/2) \sum_{j=1}^{g} \omega_j^{J} \wedge\bar{\omega}_j^{J}$.
 
When $C$ is  hyperelliptic  and otherwise general then the basic cycle
$K$ is indecomposable because
the primitive contribution of the standard regulator image of $K$
does not vanish, see \cite{3}. 
This information is found by the study  of an infinitesimal invariant of Griffiths type.
For our purpose here we need a more concrete computation on the behaviour
of the {\it real}
regulator image of $K$.
Let $ \tau := \omega_1^{J} \wedge\bar{\omega}_1^J -\omega_2^{J} \wedge\bar{\omega}_2^J$,
then   
\proclaim {2.4.1. Proposition} $ <R(K) ,\tau> \ne 0 \, .$
\endproclaim
\noindent The proof is given below in (2.6).

\subhead 2.5. Lewis' condition holds for the regulator image of ${\frak Z }$ 
\endsubhead
We take $g=3$ and $C$ general enough so that  Lieberman's result applies to $J(C)$. 
It is straightforward to check that $ \eta := \tau \wedge \omega_3^{J}$ gives
a primitive class in   $ H^{2,1}(J(C)) $ and thus $ \eta$ is
orthogonal to the coniveau $1$ filtration. 
Let  $\nu$  be the antiholomorphic form on $J \times C$
which is the pull-back of $\bar \zeta_3 = a^{\ast} \bar \omega_3^{J}$  
from the second factor $C$. 
We have

\proclaim{2.5.1. Theorem} $ <R({\frak Z }) ,\eta \wedge \nu > \ne 0 $.
\endproclaim
\demo{Proof} It is
$ <R({\frak Z }) ,\eta \wedge \nu > = $
$ <R({\frak S }) ,\eta \wedge \nu >$ -
$ <R({\frak T }) ,\eta \wedge \nu > $.
To compute we pull back $\eta \wedge \nu$  to $C \times C$ by
means of the $4$ relevant maps  $C \times C \to J \times C$.
In the straight case, that is  ${\frak S }$, the pull back of $\eta$ is the zero form,
because a $3$-form dies on the curve $W_i$, and therefore
$ <R({\frak S }) ,\eta \wedge \nu > = 0$.

The twisted embedding  with image $T_1$  is the map 
$ g : C\times C \to J \times C$ defined as  $g(x,t) = (g_1(x,t),t)$,
where  $g_1 : C\times C \to J$ is the function $g_1(x,t) = class(x+t-w_1)$. 
It is   $g_1^{\ast}(\omega_n^{J}) =\zeta_ n^{1} + \zeta_ n^{2} $, 
here the superscript $i$ indicates pull-back under the i-th projection
$C\times C \to C$. Keeping this notation we obtain:
$$ g^*(\eta \wedge \nu)  =
 (\zeta_1^{1} \wedge\bar{\zeta}_1^{1} -\zeta_2^{1}  \wedge\bar{\zeta}_2^{1}) 
\wedge (\zeta_3^{2} \wedge\bar{\zeta}_3^{2}) \,  +  
\sum_i \alpha_i^{1} \wedge \beta_i^{2}$$

\noindent where the forms $\beta_i^{2}$ are pull-back from the second factor
of forms $\beta_i$ with the property  $\int_{C} \beta_i =0$ .
Our aim is to check the non vanishing of 
$$
 \int_{T_1} 
 \eta \wedge \nu \,\log|G_1| \ +  \int_{T_2} 
 \eta \wedge \nu \,\log|G_2| \ 
$$
going to $C \times C$ 
this is twice 
 
 $$
 \int_{C \times C}  g^*(\eta \wedge \nu) \,\log|F| \ ,
$$
Now  $F$ is the pull-back of the rational function $f$ from the first factor $C$ 
and therefore we can use the product structure to calculate the integrals.
One has

$$\align
&  \int_{C \times C}  g^*(\eta \wedge \bar{\zeta_3}) \,\log|F|  =    \\
& ( \int_{C}  
 (\zeta_1^{1} \wedge\bar{\zeta}_1^{1} -\zeta_2^{1}\wedge\bar{\zeta}_2^{1}) \log|f| )(\int_{C}  
 \zeta_3^{2} \wedge\bar{\zeta}_3^{2})    + 
\sum_i  (\int_{C}  \alpha_i^{1} \log|f|)( \int_{C}  \beta_i^{2})  \, =  \\
& ( \int_{C}  
 (\zeta_1^{1} \wedge\bar{\zeta}_1^{1} -\zeta_2^{1}\wedge\bar{\zeta}_2^{1}) \log|f| )(\int_{C}  
 \zeta_3^{2} \wedge\bar{\zeta}_3^{2})  \quad  \quad  \quad \quad \quad  .   \endalign  $$

\noindent Up to a nonzero constant this  is precisely $ <R({K}) ,\tau>$,
the regulator pairing for the basic cycle,  hence  it is not  zero by proposition 2.4.1. 
\enddemo

\subhead 2.6. Useful integrals and the proof of 2.4.1
\endsubhead

 The proof of 2.4.1 depends to a large extent on  a reduction process
to the case of elliptic curves, as we explain next.
 
Let $E_{\lambda}$ be the elliptic curve with affine equation
$y^2 = x(x-1)(x-{\lambda})$.  Then we define $f_{\lambda} := x$
as a rational function on $E_{\lambda}$. Let 
$ \theta_{\lambda} = (i/2)   {\omega_{\lambda} } \wedge\bar{\omega}_{\lambda} $
be the invariant volume form on $E_{\lambda}$.
We define    $I({\lambda}) := \int_{E_{\lambda}}  \log|f_{\lambda}|
\theta_{\lambda} $.
 
\proclaim{2.6.1. Proposition} $I({\lambda})$  varies with ${\lambda}$.
\endproclaim
\demo{Proof}
Multiplication of $x$ by ${\lambda}^{-1}$ shows that 
${\lambda}$ and ${\lambda}^{-1}$  determine the same curve
$E$, say.   On $E$ the volume forms coincide,
while $f_{\lambda} =  {\lambda} f_{{\lambda}^{-1}}$. In this way 
$$I({\lambda})= \int_{E}  \log|{\lambda}| \theta +  I({\lambda}^{-1})
= \log|{\lambda}|+ I({\lambda}^{-1})$$

\noindent hence $I({\lambda})$ cannot be constant.
\enddemo 

	We move to the case of genus $2$ and define $\tau $ as in $2.4$.
Consider the hyperelliptic map $f: C \to \Bbb P^1$ and define
$I(f,\tau) := \int_{C}  \log|f| \tau $.
 
\proclaim{2.6.2. Proposition} If $C$ is general and of genus $2$ then $I(f,\tau) \ne 0$.
\endproclaim 
\demo{Proof}
We prove this for the genus two bielliptic curve $C$ which is a double cover of
$E(1):= E_{\lambda_1}$  and of $E(2):= E_{\lambda_2}$.
Consider the diagram
$$
\CD
E_2 @< k_2 << C @> k_1 >> E_1
\\
@V    f_2 VV @V    f VV @V  f_1  VV    
  \\ 
\Bbb P^1 @< h << \Bbb P^1 @> h  >> \Bbb P^1 
\endCD $$ 
\noindent  
Here $f_i$ is ramified over $(0,1,\infty, \lambda_i )$, $h$ is 
the double cover ramified over 
$ \lambda_1 $ and $ \lambda_2 $, and 
$f:  C \to \Bbb P^1$ is the hyperelliptic cover  ramified at
$ h^{-1}(\{0,1,\infty \}) $. On the range of $h$
we have already fixed a standard parameter,
we choose a standard parameter on the domain of $h$
so that $0$ maps to $0$, and similarly for $1$ and for $\infty$.
In this manner $f$ is a well defined rational function on $C$,
and we denote by $\bar f$ its transform under the involution
of $\Bbb P^1$ associated with $h$. Consider the
composition $g:=hf$, then $g$ is a rational  function 
of degree $4$ on $C$ and in fact $f \bar f = cg$, $c$  a non-zero  constant.

Given a rational function $r$ on $C$ we set for $i = 1$ or $= 2$ 
$$I(r,i,C) := \int_{C}k_i^{*}(\theta_{i})  \log|r|  \, , \quad
 \qquad  I(r, \theta_{1} -  \theta_{2} ) := I(r,1,C) - I(r,2,C) $$

\noindent Using the previous notations we see that 
$I(g,i,C) =  I({\lambda}_i)$,  
and therefore  in general it is $I(g, \theta_{1}  -  \theta_{2} ) \neq  0$.
We have
$$I(f, \theta_{1}  -  \theta_{2} )   + I(\bar f, \theta_{1}  -  \theta_{2} ) = 
I(g, \theta_{1}  -  \theta_{2} ) +   
log|c| \int_{C} ( k_1^{*}(\theta_{1})  - k_2^{*}(\theta_{2})) \, =  
I(g, \theta_{1}  -  \theta_{2} )$$

\noindent and therefore   $I(f, \theta_{1}  -  \theta_{2} ) \neq 0 $ 
for the general bielliptic
curve C. Recall now that  $ (1/2 )(k_1^{*}(\theta_{1}) - k_2^{*}(\theta_{2}))$ 
can be taken to be $\tau_C$  and then conclude 
that $ I(f,\tau)  \neq 0 $.  
\enddemo

	To complete the proof
of 2.4.1 we show now that it holds for the case of a curve $C$ which 
is an unramified double cover $ \pi : C \to G$,
where $G$ is a general curve of genus $2$. 
To begin we remark that  $\pi^{\ast} \tau_G$ is a 
multiple of a form on $C$ which can be taken  as $\tau_C$. 
There is a diagram 

$$
\CD
 C @> {\pi} >> G @< {\pi} << C
\\
 @V  f VV @V  f_G  VV  @V  k_E  VV  
  \\ 
\Bbb P^1   @> h >> \Bbb P^1  @< f_E << E
\endCD $$ 

\noindent  here  $ \pi : C \to G$ is \'etale,
 $ f_G :   G  \to  \Bbb P^1$  is the
hyperelliptic cover, ramified at $6$ points, $\{ 0,1, \infty, \lambda, a_1, a_2\}$.
The  cover $h$ is ramified at  $a_1, a_2$ and  $f_E :E \to \Bbb P^1 $
is ramified at the other points.  The map
$ f : C \to \Bbb P^1$ is also an hyperelliptic cover,
it is ramified over the pull-back via $h$ of the ramification points 
of   $f_E$, and finally   $ k_E : C \to E$ is ramified along the pull back
via $f_E$ of the ramification of $h$.
We set here $g:= hf$ and define the parameter on the domain  of
$h$ in such a way that $h(0) = 0$, $h(1) =1 $,  $h(\infty) = \infty $.
Using  2.6.2. for $G$ we conclude by a similar argument 
that 2.4.1 holds.

\heading 3.     Natural higher cycles on the general jacobian of genus 3 
 \endheading
\subhead 3.1.      The 4-configuration
\endsubhead 

By 4-configuration we mean here a certain cycle in $CH^3(J(C),1)$  which is supported
on $4$ copies of the curve $C$.
	The configuration depends on the choice of points 
$a'$ and $a''$  and $p'$ and $p''$ on $C$ with the condition that 
$class ((a'+a'')-(p'+p'')) = \epsilon$ is a torsion $2$ element in $J(C)$.
More precisely we choose a rational function $f$ with $div(f) =  2((a'+a'')-(p'+p''))$.
We embed $C$ in $Pic^3$  in $4$ ways as follows.
Let $ i(y,z): C \to Pic^3  $ be the map  $ i(y,z) (x) =  x+y+z$,
and   $C(y,z) := i(y,z) (C)$,   and let 
$ j : C \to Pic^3  $ be the map  $ j (x) =  -x+2(a'+a'')$, $G := j(C)$.
We consider $ C(a',a'')$, $C(p',p'')$ and $G$. Translation by ${\epsilon}$
maps $ C(a',a'')$ to $C(p',p'')$, and we take $G_{\epsilon}$ to be 
the image of $G$.
	We shall use  the convention that  $f$ represents the rational function
on each of the preceding curves which maps to $f$ under the chosen isomorphism with
$C$ and thus we set $ Z_1  := C(a',a'') \otimes f$,
$ Z_2  := G \otimes f $, 
$ Z_3  := C(p',p'') \otimes f$,  
$ Z_4  := G_{\epsilon} \otimes f$. 
 
\proclaim{Proposition} $Z := \sum_{i=1}^{4} (-1)^i Z_i$    is   a higher cycle. 
\endproclaim 
The only possible difficulty is to see where
the curves intersect. Now  $ C(a',a'')$ intersects  $G $ in two points,
on both curves  the points come from  $a'$ and $a''$ under the isomorphism with $C$,
but the point which comes from $a'$ in $G$ comes from $a''$ in $ C(a',a'')$
and conversely. A similar statement holds for $C(p',p'') \cap G$, and then intersections
with $G_{\epsilon}$ can be recovered by using ${\epsilon}$-symmetry.
Note that if $C$ is not hyperelliptic then
$C(a',a'')\cap  C(p',p'')  = \emptyset $ and $ G \cap G_{\epsilon} = \emptyset $.  

\subhead 3.2. The hyperelliptic configuration is a special kind of
4-configuration  \endsubhead  

	It is apparent that there is a $1$-dimensional family 
of $4$-configurations associated with a given
$2$ torsion class $\epsilon$ for a fixed curve $C$ of genus $3$. We specialize 
the curve to be hyperelliptic. We degenerate also the given pencil $g^1_4$, 
by the request that $a'$  and $p'$, say, 
move to become the same point $t$ on the hyperelliptic
curve $C$, and therefore $a''$ becomes the Weierstrass point $w_1$  and $p''$ is $w_2$.
In this manner the rational function $f$ must be the Weierstrass one 
and the pencil $g^1_4$ is 'degenerate', namely  $g^1_2$ $ + 2t$.
Strictly speaking our construction should be called deformation {\it from} the hyperelliptic
curve. In this hyperelliptic case we choose to denote by $Z(t)$ the higher 
cycle which is constructed according to the rule of the 4-configuration 
so to make explicit its dependence on $t$. It is simple to check that
$Z(t)$ is indeed the hyperelliptic configuration
$K- K_t$  of section $1.2$, and therefore it is indecomposable.  
This is best seen by identifying $ Pic^3  \to Pic^1(C)= J(C) $  
by first mapping $ Pic^3 \to Pic^{-3}$ and then going from $Pic^{-3}$ to
$Pic^{1}$ by adding the divisor $2(t+w_1)$.
In this manner   $G$ and $G_{\epsilon}$ become
$W_1$  and  $W_2$ respectively, while
the curve $ C(t,w_1) $ goes to  $C_t$ 
and  $ C(t,w_2) $  to $C_{t+ \epsilon}$.  
There is apparently a minor problem here, in that the chosen
isomorphism of $C$ with $C_t$ 
differs from the isomorphism $C \to C(t,w_1) \to C_t$,
by the hyperelliptic involution. This fact is of no
consequence, because
the Weierstrass function is of course invariant under
the involution.

\remark{3.2.1. Remark} We have seen that on the general Jacobian
the $4$-configuration 
yields $4$ curves and $8$ points of intersection.
Each curve does not intersect its $\epsilon$-translate,
and meets the remaining curves each in $2$ points.
Each point belongs to $2$ curves. On the other hand   
when the configuration becomes  the hyperelliptic one
just described, then the points of intersection become $4$ in all,
every curve meets every other  in $2$ points, but now 
each point belongs to $3$ curves and each curve to $3$
points. 
\endremark      

\bigskip
I do not know how to turn into a proof the heuristic observation  that the indecomposability
of the hyperelliptic configuration should  yield indecomposability
for the general 4-configuration, still I think that there is enough evidence for: 

\proclaim{3.2.2  Conjecture} The general 4-configuration is indecomposable.
\endproclaim

A strong motivation in this direction is  Fakhruddin's work \cite {4},
indeed our theorem (2.3.1) and the conjecture may both be 
interpreted as a kind of extension to higher Chow groups of some of Fakhruddin's results. 

\heading 4. Translations act non trivially on $CH_{\text{ind}}^{g+1}(J(C),2)$. \endheading

 	We consider now
the second higher Chow group $ CH^{p}(X,2) \simeq H^{p-2}(X,\Cal K_{p})$,
according to the terminology of \cite {8} the group
of decomposable cycles is here
$$
CH_{\text{dec}}^k(X,2) := Im \big\{ K_2(\Bbb C) {\otimes}  CH^{k-2}(X)
\longrightarrow CH^k(X,2) \big\} \quad,
$$
while the {{\it} indecomposable\/} group is
$$
CH_{\text{ind}}^k(X,2;{\Bbb Q}) := 
\big( CH^{k}(X,2)\big/ CH_{\text{dec}}^k(X,2)\big) {\otimes}
{\Bbb Q} \quad.
$$ 

It is known that translations on an elliptic curve $E$ act
trivially on  $CH_{\text{ind}}^2(E,2)$, see \cite {5,  3.10}. 
We show that on the contrary translations on a genus $2$ bielliptic jacobian $J(G)$
operate  non trivially on $CH_{\text{ind}}^3(J(G),2)$.
Our procedure is similar to the one used above for the first higher Chow group.   
The hyperelliptic configuration is replaced
by a certain B-configuration, which is then shown to be indecomposable
by checking Lewis' condition on a cycle  ${\frak B }$ of $CH^{3}(J(G)\times G,2)$.
 
\subhead 4.1.  The bielliptic  configuration
\endsubhead

	Spencer Bloch  defined and studied for the first time higher regulator maps. 
In his seminal memory \cite {2} he constructed certain
elements $S_b \in \Gamma(E,\Cal K_{2}))$
associated with a point $b$ of finite order on an elliptic curve $E$.
It is one consequence of his deep work 
that the real regulator image of $S_b$ is not trivial for some curves with complex multiplication,
and therefore that it is not trivial in general. In \cite {3} it is shown that moreover  $CH_{\text{ind}}^2(E,2)$
is not finitely generated.

Consider a bielliptic curve $G$ of genus $2$ with 
associated map $\delta_G :G \to E_1 $, and let
$a: G \to J(G)$ be the Abel Jacobi map.
In this way  $Z(b) :=  a_{\ast} \delta_G^{\ast}  (S_b)$ is a cycle in
$CH^{3}(J(G),2)$.
Translation by an element $t \in Pic^0 (G)$  maps  $Z(b)$ to the cycle $Z_t(b)$, our aim
is to prove 
\proclaim{4.1.1. Theorem}
The bielliptic configuration $B(t) := Z_t(b) - Z(b)$ is indecomposable in
general.
\endproclaim

Note that $B(t)$ has trivial regulator image.

\demo{Proof}  Consider the cycle $G \times Z(b)$ as an element
in the higher Chow group of  $G \times G $.  
The straight embedding $\sigma := id \times a:  G \times G \to   G \times J(G)$
maps it to   ${\frak S } := \sigma_{\ast} (G\times Z(b))$ 
in $CH^{3}(G \times J(G),2)$.
The twisted  embedding   $\tau (t,x) := (t,a(x)+(t-w_1)))$  gives instead
${\frak T }:= \tau_{\ast} (G\times Z(b)) $, with section
${\frak T }_{\ast}(t) = Z_t(b)$, and therefore $B(t)$ is
the section at $t \in G$ of ${\frak B }:=  {\frak T } -    {\frak S }$.

	We use the same type of notations as we did in part $2$, in particular  the holomorphic form
${\omega}_i^{J}$ comes from $E_i$, we need to consider also the forms  
$\nu :=\bar{\omega}_1^{J}\wedge{\omega}_2^{J} $ 
on $J(G)$  and $\bar{\zeta_2}$ on $G$.
The  procedure of  $2.5.1$ gives here again  
$<R({\frak B }),\bar{\zeta_2} \wedge \nu > \neq 0$. 
The Neron Severi space of divisors with rational coefficients on $J(G)$
is isomorphic to the same space on the product of
the two associated elliptic curves.     
On the general bielliptic jacobian
$\nu$ is orthogonal to the Neron Severi group, because it is orthogonal
to the elliptic curves. For this reason
the K\"unneth projected image of $R({\frak B })$ yields
a non trivial element in  $H^{\{3,2,2\}}(J(G)) \big/ H^{\{1,0,0\}} (J(G))$.
A look at the proof of the main theorem of \cite {8} shows that
this fact implies that the general section ${\frak B }_{\ast}(t)$ is indeed indecomposable.
\enddemo

\enddocument